\documentclass[twocolumn,11pt]{article}
\usepackage{times}
\newtheorem{definition}{Definition}

\newtheorem{theorem}[definition]{Theorem}

\newcommand{\eop}{\hfill $\sqcap\!\!\!\!\sqcup$} 
\begin{document}
\global\def\refname{{\normalsize \it References:}}
\baselineskip 12.5pt
%
%
%
\title{\LARGE \bf Finite velocity planar random motions driven by inhomogeneous fractional Poisson distributions}

\date{}

\author{\hspace*{-10pt}
\begin{minipage}[t]{2.7in} \normalsize \baselineskip 12.5pt
\centerline{ROBERTO GARRA}
\centerline{Sapienza University of Rome}
\centerline{Dipartimento di Scienze di Base}
\centerline{ed Applicate per l'Ingegneria}
\centerline{Via A. Scarpa, Rome}
\centerline{ITALY}
\centerline{roberto.garra@sbai.uniroma1.it}
\end{minipage} \kern 0in
\begin{minipage}[t]{2.7in} \normalsize \baselineskip 12.5pt
\centerline{ENZO ORSINGHER}
\centerline{Sapienza University of Rome}
\centerline{Dipartimento di Scienze Statistiche}
\centerline{Piazza Aldo Moro, Rome}
\centerline{ITALY}
\centerline{enzo.orsingher@uniroma1.it}
\end{minipage}
%
%
\\ \\ \hspace*{-10pt}
\begin{minipage}[b]{6.9in} \normalsize
\baselineskip 12.5pt {\it Abstract:}
In this paper we study finite velocity planar random motions with an infinite number of possible directions, where
the number of changes of direction is randomized by means of an inhomogeneous fractional Poisson distribution. We first discuss
the properties of the distributions of the generalized fractional inhomogeneous Poisson process. Then we show
that the explicit probability law of the planar random motions where the number of changes of direction is governed by this fractional
distribution can be obtained in terms of Mittag-Leffler functions. We also consider planar random motions with random velocities obtained from
the projection of random flights with Dirichlet displacements onto the plane, randomizing the number of changes of direction with a suitable
adaptation of the fractional Poisson distribution.
\\ [4mm] {\it Key--Words:}
Fractional Poisson processes, Finite velocity random motions, Mittag-Leffler functions
\end{minipage}
\vspace{-10pt}}

\maketitle

\thispagestyle{empty} \pagestyle{empty}
%
%
\section{Introduction}
\label{S1} \vspace{-4pt}
Planar random motions at finite velocity have been studied
by many authors, with different approaches.
Due to the difficulty of considering a continuous spectrum of infinite possible directions
of random motions on the plane, many studies have been devoted to the analysis
of cyclical planar random motions with a finite number of possible directions (see for example \cite{dic} and \cite{enzo}).
Planar random motion at finite velocity with an infinite number of uniformly distributed orientations of
displacements has been studied by several researchers over the years. In particular Orsigher and Kolesnik in \cite{kol} studied a planar random
motion with an infinite number of directions whose probability law is governed by the
damped wave equation. In this work the number of changes of directions
is given by an homogeneous Poisson process.

In the framework of fractional point processes, a number of papers have been devoted to the analysis of different fractional Poisson-type processes (see for example \cite{Beghin, Laskin, g4}). An interesting application of the fractional Poisson distribution (introduced in \cite{Beghin})
\begin{equation}
P\{N_{\alpha}(t)= k\}= \frac{1}{E_{\alpha,1}(\lambda t)}\frac{(\lambda t)^k}{\Gamma(\alpha k+1)},
\end{equation}
was discussed in recent papers about random flights with Dirichlet distributed displacements (see \cite{ale3, Ale, SPA}). In these papers the fractional Poisson distributions play a key-role in order to find in an explicit form the unconditional probability law of random flights in dimension $d>2$ and the corresponding governing partial differential equations.

In this paper we introduce a generalized fractional Poisson distribution with time-dependent rate and we discuss some of its main properties. Then, we study planar random motions with an infinite number of possible directions, where the number of changes of directions is given by the inhomogeneous fractional Poisson process.
We will show that explicit probability law can be expressed in terms of Mittag-Leffler functions. Finally we consider planar random motions with finite
velocity obtained by the projection of random flights with Dirichlet displacements studied in \cite{Ale} onto the plane.

\subsection{Non-homogeneous fractional Poisson distributions}
\vspace{-4pt}
Here we introduce the following generalization of the fractional Poisson distribution, given by
\begin{equation}\label{uno}
P\{N_{\alpha}(t)= n\}= \frac{1}{E_{\alpha,1}(\Lambda(t))}\frac{\left(\Lambda(t)\right)^n}{\Gamma(\alpha n+1)},\quad n \geq 0
\end{equation}
where $\alpha\in (0,1]$,
\begin{equation}
E_{\alpha,1}(t)= \sum_{k=0}^{\infty}\frac{t^k}{\Gamma(\alpha k+1)},
\end{equation}
is the classical one-parameter Mittag-Leffler function and
\begin{equation}
\Lambda(t)= \int_0^t\lambda(s)ds.
\end{equation}
It is immediate to check that, for $\alpha = 1$, the distribution
(\ref{uno})
coincides with the distribution of a non-homogeneous Poisson process with time-dependent
rate $\lambda(t)$. Moreover, for $\lambda(t)= \lambda= const.$ and $\alpha \in (0,1)$, this distribution coincides with one of the
alternative fractional Poisson processes discussed by Beghin and Orsingher in \cite{Beghin}.

We now discuss the meaning of this generalization of the Poisson process in the theory of counting processes.

First of all, we should remark the relation of distribution
(\ref{uno}) with weighted Poisson processes.
We recall that the probability mass function of a weighted
Poisson process is of the form (see for example \cite{Bala}
and references therein)
\begin{equation}
P\{N^w(t)=n\}= \frac{w(n)p(n)}{E[w(N)]},\quad n\geq 0,
\end{equation}
where $N$ is a radom variable with a Poisson distribution $p(n)$,
$w(\cdot)$ is a non-negative weight function with non-zero, finite
expectation, i.e.
\begin{equation}
0<E[w(N)]= \sum_{n=0}^{\infty}w(n) p(n) <\infty.
\end{equation}
Then, we recognize that the distribution (\ref{uno}) is a weighted
Poisson process with time-dependent rate $\lambda(t)$ and weights
$w(n)= n!/\Gamma(\alpha n+1)$.

We observe that the probability generating function of the r.v. with distribution (\ref{uno}), reads
\begin{equation}
G(u;t)= \sum_{n=0}^{\infty}u^n P\{N_{\alpha}(t)= n\}=
\frac{E_{\alpha,1}(u\Lambda(t))}{E_{\alpha,1}(\Lambda(t))}
\end{equation}
and satisfies the following fractional differential equation
\begin{equation}
\frac{d^{\alpha}}{du^{\alpha}}G(u^{\alpha};t)=\Lambda(t)G(u^{\alpha};t),
\end{equation}
where
\begin{equation}
\frac{d^{\alpha}f}{du^{\alpha}}=
 \frac{1}{\Gamma(1-\alpha)}\int_0^{u}(u-s)^{-\alpha}
 \frac{d}{ds}f(s) \, ds, \, u>0,
 \end{equation}
is the Caputo fractional derivative of order $\alpha \in(0,1)$ (see
e.g. \cite{Pod}).

 We finally recall that in \cite{noi}, the authors studied a state-dependent fractional Poisson process, namely $\widehat{N}(t)$, $t \ge 0$.
        In view of the previous analysis, the non-homogeneous counterpart of this process
        has univariate probabilities given by
        \begin{equation}
                    \label{p2}
                    \Pr\{\widehat{N}(t)= j\}=\frac{\frac{(\Lambda(t))^j}{\Gamma(\alpha_j j+1)}
                    \frac{1}{E_{\alpha_j,1}(\Lambda(t))}}{\sum_{j=0}^{+\infty} \frac{(\Lambda(t))^j}{\Gamma(\alpha_j j+1)}
                    \frac{1}{E_{\alpha_j,1}(\Lambda(t))}}, 
                \end{equation}
        where $0<\alpha_j\leq 1$, for all $j \geq 0$. This distribution leads to an interesting form of counting process, where the fractional
        weights depend on the state. The statistical applications of this approach should be object of further research.

\section{Planar random motions with finite velocity and infinite directions: main results}
In this section we recall some results obtained by Kolesnik and Orsingher in \cite{kol} about planar random motion with finite velocity with an infinite number of directions.
In their model, the motion is described by a particle taking directions $\theta_j$, $j=1, 2, \dots$, uniformly
distributed in $(0,2\pi]$ at Poisson paced times. The orientations $\theta_j$ and the governing Poisson process
$\mathcal{N}(t)$, $t\geq 0$, are assumed to be independent.
The particle starts its motion from the origin of the plane at time $t=0$
        and moves with constant velocity $c$. It changes direction
        at random instants according to a Poisson process.
        At these instants, the particle instantaneously takes a new direction $\theta$ with uniform distribution
        in $[0,2\pi)$ independently of its previous deviation.
        Therefore, after $N(t)=n$ changes of direction, the position $(X(t),Y(t))$ of the particle in the plane
        is given by
        \begin{eqnarray}
        \label{plc}
        X(t) = c\displaystyle \sum_{j=1}^{n+1}(s_j-s_{j-1})\cos \theta_j,\\
        Y(t)=  c\displaystyle \sum_{j=1}^{n+1}(s_j-s_{j-1})\sin \theta_j,\label{prc}
        \end{eqnarray}
        where $\theta_j$, $j= 1, \dots,n+1$, are independent random variables
        uniformly distributed in $[0,2\pi)$, $s_j$ are the instants
        at which Poisson event occurs, $s_0=0$ and $s_{n+1}= t$. By means of (\ref{plc}) and (\ref{prc}),
        the conditional characteristic function of the random vector $(X(t),Y(t))$ can be written as follows
        \begin{eqnarray}\label{plc1}
            \nonumber &E\{e^{i\alpha X(t)+i\beta Y(t)}|N(t)=n\}\\
            &=\frac{2^{n/2}\Gamma(\frac{n}{2}+1)}{\left(ct\sqrt{\alpha^2+\beta^2}\right)^{n/2}}
            J_{\frac{n}{2}}\left(ct\sqrt{\alpha^2+\beta^2}\right), \\
            \nonumber &n\geq 1, \: (\alpha,\beta) \in R^2,
        \end{eqnarray}
        Then by inverting (\ref{plc1}), the following conditional distribution can be found (see formula (11) of \cite{kol})
        \begin{eqnarray}
            \label{pt1}
             &P\{X(t)\in dx, Y(t)\in dy|N(t)=n \}\\
            \nonumber &=\frac{n \left(c^2t^2-x^2-y^2\right)^{\frac{n}{2}-1}}{2\pi(ct)^n} dx \, dy,
        \end{eqnarray}
        for $x^2+y^2<c^2t^2$, $n\geq 1$.
        In the model of Kolesnik and Orsingher, the changes
        of direction are driven by an
        homogeneous Poisson process, such that the absolutely continuous
        component of the unconditional distribution of $(X(t),Y(t))$ reads
        \begin{eqnarray}
            \nonumber &P\{X(t)\in dx, Y(t)\in dy\}\\
            \nonumber &=\frac{\lambda}{2\pi c}
            \frac{e^{-\lambda t+\frac{\lambda}{c}\sqrt{c^2t^2-x^2-y^2}}}{\sqrt{c^2t^2-x^2-y^2}} dx \, dy,\label{pt}
        \end{eqnarray}
        for $x^2+y^2<c^2t^2$.

        The singular component of $(X(t),Y(t))$ pertains to the probability of no
         Poisson events. It is uniformly distributed on the circumference of radius $ct$
        and has weight $e^{-\lambda t}$.
        It has been proven that the density in (\ref{pt}) is a solution to the planar telegraph equation
        (also equation of damped waves)
        \begin{equation}
            \frac{\partial^2 p}{\partial t^2}+2\lambda\frac{\partial p}{\partial t}
            = c^2 \left\{\frac{\partial^2}{\partial x^2}+\frac{\partial^2}{\partial y^2}\right\}p.
        \end{equation}

\section{Planar motions where the number of changes of direction is given by fractional Poisson distributions}
We now apply the family of fractional-type Poisson distributions (\ref{uno})
in order to obtain explicit probability laws of planar random motions
with infinite changes of direction. In this case we randomize the number of changes of direction by means of the general family of distributions
(\ref{uno}), depending on the real parameter $\alpha\in(0,1)$ and
the particular choice of the time-dependent rate $\lambda(t)$.
This means that we are able to construct a general family of
planar random motions that includes as a special case the one studied in
\cite{kol}. Moreover as we will see, their explicit probability law can be expressed in terms of Mittag-Leffler function, suggesting the possible relation with fractional hyperbolic partial differential equation.

   \begin{theorem}\label{new}
   The probability law
   \begin{equation}
   p(x,y,t)=\frac{P\{X(t)\in dx, Y(t)\in dy\}}{dx\, dy}
   \end{equation}
   of the random vector $(X(t),Y(t))$ in (\ref{plc})-(\ref{prc}), when the number of changes of
   direction is given by the distribution (\ref{uno})
   has the following form
   \begin{eqnarray}
    \nonumber & p(x,y,t)=
   p_{ac}(x,y,t)I_{C_{ct}}+p_s(x,y,t)I_{\partial C_{ct}}\\
   &= \frac{1}{E_{\alpha,1}(\Lambda(t))}\frac{\Lambda(t)E_{\alpha,\alpha}\left(\frac{\Lambda(t)}{ct}
      \sqrt{c^2t^2-(x^2+y^2)}\right)}{2\pi \alpha ct\sqrt{c^2t^2-(x^2+y^2)}}I_{C_{ct}} \nonumber\\
      &+\frac{1}{E_{\alpha,1}(\Lambda(t))}I_{\partial C_{ct}},\quad \alpha \in(0,1]\label{bello}
   \end{eqnarray}
   where
   \begin{equation}
            C_{ct}=\{(x,y)\in R^2: x^2+y^2< c^2t^2\},
   \end{equation}
    $I_{C_ct}$ and $I_{\partial C_{ct}}$ are the indicator functions of the circle $C_{ct}$ and its boundary.
     \end{theorem}
   {\bf Proof:} The result stems from the fact that
   \begin{eqnarray}
    \nonumber p_{ac}(x,y,t)= \sum_{n=0}^{\infty}&
    P\{X(t)\in dx, Y(t)\in dy| N_{\alpha}(t)= n\}\\
    &\times P\{N_{\alpha}(t)=n \},
   \end{eqnarray}
   using the conditional distribution (\ref{pt1}) and the 
   non-homogeneous fractional Poisson distribution (\ref{uno})
   with an arbitrary time-dependent rate $\lambda(t)$.\\
   The singular component coincides with the probability of no changes
   of directions according to the distribution (\ref{uno}) and is concentrated on the boundary $\partial C_{ct}$. On the other hand it is simple to check that
   \begin{equation}
    \int\int_{C_{ct}} p_{ac}(x,y,t)dxdy= 1- \frac{1}{E_{\alpha,1}(\Lambda(t))}.
   \end{equation}

   \eop

    An interesting example is given by the special case
    in which $\Lambda(t)= \lambda t$. In this case the absolutely continuous component of the distribution has the following
    simple form
    \begin{equation}\label{due}
    p_{ac}(x,y,t)= \frac{1}{E_{\alpha,1}(\lambda t)}
    \frac{\lambda E_{\alpha,1}\left(\frac{\lambda}{c}
          \sqrt{c^2t^2-(x^2+y^2)}\right)}{2\pi c\sqrt{c^2t^2-(x^2+y^2)}},
    \end{equation}
    that, for $\alpha = 1$ coincides with the distribution of planar
    random motions discussed in \cite{kol}.

    The projection of the planar random motion with distribution (\ref{bello}) in the line can be calculated as follows
    \begin{eqnarray}
         \label{ciao} &p(x,t)=\int_{-\sqrt{c^2t^2-x^2}}^{\sqrt{c^2t^2-x^2}}p(x,y,t)dy\\
         \nonumber &=\frac{1}{\sqrt{\pi}E_{\alpha,1}(\Lambda(t))}\sum_{k=0}^{\infty}\left(\frac{\Lambda(t)}{ct}\right)^k
         \frac{\Gamma(\frac{k}{2}+1)}{\Gamma(\frac{k+1}{2})}\frac{(\sqrt{c^2t^2-x^2})^{k-1}}{\Gamma(\alpha k+1)}.
     \end{eqnarray}
     We remark that the zeroth-term of the series in (\ref{ciao}) pertains to the projection of the singular component of the distribution
     (\ref{bello}).

    The function appearing in (\ref{ciao}) can be expressed in terms of the generalized Wright function (see \cite{Kil} and
    references therein)

    \begin{eqnarray}
    \nonumber &_p\psi_q\bigg[t\bigg|
     \begin{array}{cc}
      (a_1,A_1),\dots,(a_p, A_p)\\
     (b_1,B_1),\dots, (b_q, B_q)
     \end{array}
     \bigg]\\
     \nonumber &=\displaystyle\sum_{k=0}^{\infty}\frac{\prod_{j=1}^p\Gamma(a_j+A_j k)}{\prod_{j=1}^q\Gamma(b_j+B_j k)}\frac{t^k}{k!},
    \end{eqnarray}
    where $a_j,b_j \in R$ and $A_j, B_j>0$, for all $j\in N$.
    Hence, we can write the distribution (\ref{ciao}) in the more compact form
    \begin{equation}
    \nonumber p(x,t)= \frac{_2\psi_2\bigg[\frac{\Lambda(t)}{ct}\sqrt{c^2t^2-x^2}\bigg|
         \begin{array}{cc}
          (1,1),(1, \frac{1}{2})\\
         (\frac{1}{2},\frac{1}{2}), (1, \alpha)
         \end{array}
         \bigg]}{\sqrt{\pi}E_{\alpha,1}(\Lambda(t))\sqrt{c^2t^2-x^2}}
    \end{equation}

    For $\alpha = 1$, from (\ref{ciao}), we have that
     \begin{eqnarray}
    \nonumber &p(x,t)= e^{-\Lambda(t)}\displaystyle\sum_{k=0}^{\infty}\left(\frac{\Lambda(t)}{2ct}\right)^k\frac{(\sqrt{c^2t^2-x^2})^{k-1}}{[\Gamma(\frac{k+1}{2})]^2} \\
    \nonumber &=\displaystyle\sum_{k=0}^{\infty} \frac{(\Lambda(t))^k}{k!}e^{-\Lambda(t)}\frac{k!}{[\Gamma(\frac{k+1}{2})]^2}\frac{(\sqrt{c^2t^2-x^2})^{k-1}}{(2ct)^k}\\
    \nonumber &= \displaystyle\sum_{k=0}^{\infty}P\{N_1(t)=k\}P\{X(t)\in dx| N_1(t)=k\},
    \end{eqnarray}
    where $N_1(t)$ is the non-homogeneous Poisson distribution (\ref{uno}) with $\alpha = 1$. In such a way, we can infer the conditional
    distribution of the motion performed by the projection on the line of the random planar motion previously considered.

    Moreover, for $\alpha = 1$ and $\Lambda(t)= \lambda t$, we have that
    \begin{equation}\label{sonin}
    p(x,t)= e^{-\lambda t}\sum_{k=0}^{\infty}\left(\frac{\lambda}{2c}\right)^k
             \frac{(\sqrt{c^2t^2-x^2})^{k-1}}{[\Gamma(\frac{k+1}{2})]^2},
    \end{equation}
    that coincides with formula (1.3) of \cite{ale2} where the projection of planar random motions with and infinite number of directions and homogeneous Poisson driven changes of direction was studied. The function appearing in (\ref{sonin}) is known as Sonine function.
    The probability law (\ref{ciao}) describes a finite velocity motion on the line with random velocity. Indeed, in this case the distribution is completly concentrated in $|x|<ct$ because of the projection of the singular component of the planar distribution.

    We finally consider the relation between fractional differential equations and the asbsolutely continuos component of the fractional telegraph process (\ref{due}). We observe that the function
    \begin{eqnarray}
    &f(x,y,t)= p_{ac}(x,y,t)E_{\alpha,1}(\lambda t)\\
    \nonumber &  =\frac{\lambda E_{\alpha,1}\left(\frac{\lambda}{c}
              \sqrt{c^2t^2-(x^2+y^2)}\right)}{2\pi c\sqrt{c^2t^2-(x^2+y^2)}},
    \end{eqnarray}
    can be written in terms of the variable $w= \sqrt{c^2t^2-(x^2+y^2)}$.
    Then, we have that $f(w)$ satisfies the following fractional differential equation
    \begin{equation}
    \frac{d^{\alpha}}{dw^{\alpha}}\left(w^{\alpha}f(w^{\alpha})\right)=
    \frac{\lambda}{c}\left(w^{\alpha}f(w^{\alpha})\right),
    \end{equation}
    where the fractional derivative of order $\alpha \in(0,1]$ is in the Caputo sense. Indeed, we have that
    \begin{equation}
    w^{\alpha}f(w^{\alpha})= E_{\alpha,1}\left(\frac{\lambda}{c}w^{\alpha}\right),
    \end{equation}
    that is a well--known eigenfunction of the Caputo fractional derivative.

   \section{Planar motions with random velocities where the number of changes of direction is given by fractional Poisson
   distributions}
   In this section we consider a planar random motion with random velocities. Our construction is based on the marginal distributions
   of the projection of random flights with Dirichlet displacements in $R^d$ onto $R^2$ considered by De Gregorio and Orsingher in \cite{Ale}.
   In this paper the authors considered random motions at finite velocity in $R^d$, with infinite possible directions uniformly distributed
   on the hypersphere of unitary radius and changing directions at Poisson paced times. Two different Dirichlet distributions of the displacements
   were considered and the explicit form of marginal distributions of the random flights were found.

   Here we consider planar motions with random velocities, obtained from the projection of the random flights $\mathbf{X}_{d}(t)$ and $\mathbf{Y}_{d}(t)$, $t>0$, studied in \cite{Ale}
   onto the plane. We recall from \cite{Ale} (Theorem 4, pag.695) that the marginal distributions of the projections of the processes
     onto
           $R^2$ are given by
           \begin{eqnarray}
           \nonumber &f^d_{\mathbf{X}_2}(\mathbf{x}_2,t;n)= \frac{\Gamma\left(\frac{n+1}{2}(d-1)+
           \frac{1}{2}\right)}{\Gamma\left(\frac{(n+1)}{2}(d-1)-\frac{1}{2}\right)}
           \frac{(c^2t^2-\|\mathbf{x}_2\|^2)^{\frac{n+1}{2}(d-1)-\frac{3}{2}}}{\pi(ct)^{(n+1)(d-1)-1}},\\
           \nonumber &f^d_{\mathbf{Y}_2}(\mathbf{y}_2,t;n)= \frac{\Gamma\left((n+1)(\frac{d}{2}-1)
           +1\right)}{\Gamma\left((n+1)(\frac{d}{2}-1)\right)}
           \frac{(c^2t^2-\|\mathbf{y}_2\|^2)^{(n+1)(\frac{d}{2}-1)-1}}{\pi(ct)^{2(n+1)(\frac{d}{2}-1)}},
           \end{eqnarray}
           with $\|\mathbf{x}_2\|<ct$ and $\|\mathbf{y}_2\|<ct$. Here we construct exact probability distributions of planar motions
           with random velocities, randomizing the number of changes of direction with the distribution (\ref{uno}), with a suitable choice
           of the parameter $\alpha$ depending on the dimension $d$ of the original space.

           Let us consider in detail, for the sake of clarity, the planar motion obtained by the projection of the random flight $\mathbf{Y}_{d}(t)$
           of \cite{Ale}.
           We randomize the number of changes of directions by means of the following adaptation of the distribution (\ref{uno})
        \begin{equation}
         P\{N_{d}(t)= n\}= \frac{1}{E_{\frac{d}{2}-1,\frac{d}{2}}(\Lambda(t))}\frac{\left(\Lambda(t)\right)^n}{\Gamma((n+1)(\frac{d}{2}-1)+1)},\nonumber
        \end{equation}
        where $n\geq 0$ and $d$ is the dimension of the original space of the random flight that we are projecting onto the plane.
        Then, we have that the unconditional probability law of the planar motion can be obtained as follows
        \begin{eqnarray}
        \nonumber &p(\mathbf{y}_2,t)=  \displaystyle \sum_{n=0}^{\infty}
            P\{Y_1(t)\in dy_1, Y_2(t)\in dy_2| N_{d}(t)= n\}\\
            \nonumber &\times P\{N_{d}(t)=n \}\\
            \nonumber & = \frac{(c^2t^2-\|\mathbf{y}_2\|^2)^{\frac{d}{2}-2}}{\pi (ct)^{d-2}} \frac{E_{\frac{d}{2}-1,\frac{d}{2}-1}\left(\frac{\Lambda(t)}{(ct)^{d-2}}(c^2t^2-\|\mathbf{y}_2\|^2)^{\frac{d}{2}-1}\right)}{E_{\frac{d}{2}-1,\frac{d}{2}}(\Lambda(t))}.
        \end{eqnarray}
        An interesting case is for $d= 4$, where we have that
        \begin{equation}\label{bell}
        p(\mathbf{y}_2,t) = \frac{\Lambda(t)}{\pi(ct)^2}\frac{exp\bigg\{\frac{\Lambda(t)}{c^2t^2}(c^2t^2-\|\mathbf{y}_2\|^2)\bigg\}}{e^{\Lambda(t)}-1}.
        \end{equation}
        In equation (\ref{bell}) we used the fact that
        \begin{eqnarray}
        \nonumber & E_{1,2}(t)= \frac{e^t-1}{t}\\
        \nonumber & E_{1,1}(t)= e^t.
        \end{eqnarray}

        A similar derivation of the explicit probability law can be given for the other case considered in \cite{Ale}.

\end{document}